# MULTI-DIMENSIONAL PROPERTIES OF ONE-DIMENSIONAL DISCRETE FOURIER TRANSFORM

by Andrew V. BATRAC

E-mail:   ndbatrac@hotmail.com

polartechnics@polartechnics.com.au

Company: Polartechnics Ltd, Sydney

www.polartechnics.com.au

**Abstract.** A mathematical relation between elements of one- and multi-dimensional discrete Fourier transforms (DFT) is found. A method of analysing the multi-dimensional data by their single one-dimensional (1-D) DFT is offered. An experiment of filtering a two-dimensional image using a single 1-D DFT is carried out.

Subject terms: Fourier transforms, image processing.

# 1. Introduction

It is widely accepted that the multi-dimensional (M-D) discrete Fourier transforms (DFT) require an N number of calculations of one-dimensional (1-D) DFT along each of the dimensions, where N depends on both the number of dimensions (axes) and the length of data along each of the axes. Thus, for a square image comprising of *(N x N)* pixels it is required to calculate *2N* 1-D DFT [Refs. 1-4]. Actual implementation usually involves calculating either *2N* 1-D fast Fourier transforms (FFT) or a custom-made 2-D FFT, which generally does the same but *in place* [Ref. 5].

The present work has resulted from an attempt to avoid the use of 2-D FFT for the image filtering. The philosophy behind it was as follows: take *e.g.* the computer memory. It is linear and one-dimensional. We do not need to create a special M-D memory to store the M-D arrays. Indeed they are stored column after column, row after row, page after page etc., in a single memory line. Therefore, we can say that an element of a 1-D line of (say spatial) data may serve at the same time as an element of M-D (spatial) array. Believing the spatial and frequency domains to be symmetrical, it is reasonable to expect similar behaviour of a line of 1-D frequency data.

# 2. Assumptions

Any continuous function as analysed further below, be it 1-D or M-D, is considered to be regular, periodic, band-limited, over-sampled and therefore fully recoverable. Let $Euler[\theta] \equiv \exp[-2\pi \cdot \theta \cdot \sqrt{-1}]$, so that DFT pairs may be written as:

1-D for N samples:
$$F(u) \equiv \sum_{x=0}^{N-1} Euler[\frac{xu}{N}] f(x) \qquad (1.1)$$

$$f(x) = \sum_{u=0}^{N-1} Euler[-\frac{xu}{N}] \frac{F(u)}{N} \qquad (1.2)$$

2-D for *(N x M)* samples:
$$F(u,v) \equiv \sum_{y=0}^{M-1}\sum_{x=0}^{N-1} Euler[\frac{xu}{N} + \frac{yv}{M}] f(x,y) \quad (1.3)$$

$$f(x,y) = \sum_{v=0}^{M-1}\sum_{u=0}^{N-1} Euler[-\frac{xu}{N} - \frac{yv}{M}]\frac{F(u,v)}{NM} \quad (1.4)$$

Each of these DFTs obeys the translation (shifting) rule, so, *e.g.* Eq. (1.3) may be

re-written as:
$$F(u-\alpha, v) \equiv \sum_{y=0}^{M-1}\sum_{x=0}^{N-1} Euler[\frac{x}{N}(u-\alpha) + \frac{yv}{M}] f(x,y) \quad (1.5)$$

where $\alpha$ is any real number. (The continuous function recovery (restoration) from its samples does not depend on the position of the first sample, so $\alpha$ may be any real number.)

## 3. Dissertations

## 3.1 Two-dimensional Case

Let a continuous function $f(x_c, y_c)$ be properly sampled by an *(N x M)* grid to form an array of samples *f(x,y)*, where $x = 0, 1, 2, \ldots N-1$, and $y = 0, 1, 2, \ldots M-1$. Let *N* be the number of columns, and *M* the rows. Figure 1, with a (4 x 3) array, illustrates the following steps. The product *(N x M)* is limited, so we can always re-number the 2-D array *f(x,y)* as a 1-D array *f(k)*, where
$$k = Ny + x \quad (2.1)$$

For any given $k = 0, 1, 2, \ldots NM-1$, there is a unique pair of *(x,y)* and vice versa. Summation across the array may be done either as a 1-D sum, or as 2-D double sum:

$$\sum_{k=0}^{NM-1}(\ldots) = \sum_{y=0}^{M-1}\sum_{x=0}^{N-1}(\ldots) \quad (2.2)$$

Let us take a 1-D DFT (1.1) of the *f(k)*:
$$F(j) = \sum_{k=0}^{NM-1} Euler[\frac{jk}{NM}] f(k) \quad (2.3)$$

This produces a 1-D array $F(j)$, where $j = 0, 1, 2, \ldots NM-1$. It again may be re-arranged (re-numbered) as 2-D array of same size as $f(x,y)$, so $j$ becomes either $(Nv + u)$, which leads to nothing special, or $(Mu + v)$, which does the trick:

$$j = Mu + v \qquad (2.4)$$

Substituting Eqs. (2.1), (2.2) and (2.4) turns the Eq. (2.3) into:

$$F(j) = \sum_{y=0}^{M-1}\sum_{x=0}^{N-1} Euler[\frac{(MNyu + Mxu + Nyv + xv)}{NM}]f(x,y) \qquad (2.5)$$

and further:

$$F(j) = \sum_{y=0}^{M-1}\sum_{x=0}^{N-1} Euler[yu]Euler[\frac{x}{N}(u+\frac{v}{M})+\frac{yv}{M}]f(x,y) \qquad (2.6)$$

Having $u$ and $y$ as integers, $Euler[yu] = 1$ (always), and therefore may be omitted:

$$F(j) = \sum_{y=0}^{M-1}\sum_{x=0}^{N-1} Euler[\frac{x}{N}(u+\frac{v}{M})+\frac{yv}{M}]f(x,y) \qquad (2.7)$$

The left-hand part of Eq. (2.7) is still the 1-D DFT of $f(x,y)$. The right-hand part is identical to Eq. (1.5) with $\alpha = -v/M$. Finally, for $j = Mu + v$:

$$F_{1-D}(Mu + v) = F_{2-D}(u + v/M, v) \qquad (2.8)$$

## 2.2 Multi-Dimensional Case

Let $N_1, N_2, \ldots, N_m$ be the dimensions of an M-D array of samples $f(x_1, x_2, \ldots, x_m)$, where $x_i$ is the index along the $i$-th axis. For convenience let $N_0 = N_{m+1} = 1$, and the size of the array be:

$$Q = \prod_{j=0}^{m+1} N_j \qquad (3.1)$$

Repeating all the steps as before: $f(x_1, x_2, \ldots, x_m) \to f(z)$, where

$$z = \sum_{i=1}^{m} x_i \prod_{j=0}^{i-1} N_j \qquad (3.2)$$

Applying 1-D DFT:

$$F(w) = \sum_{z=0}^{Q-1} Euler[\frac{zw}{Q}]f(z) \qquad (3.3)$$

Repeating the trick, *i.e.* making the fastest changing index of $f(x_1, x_2, ..., x_m)$ to be the slowest changing index of a new array, and so forth:

$$\text{Let} \quad w = \sum_{k=1}^{m} u_k \prod_{h=k+1}^{m+1} N_h \quad (3.4)$$

$$\text{and} \quad \sum_{k=0}^{Q-1}(.....) = \sum_{Xm=0}^{Nm-1} .... \sum_{X2=0}^{N2-1} \sum_{X1=0}^{N1-1}(.....) \quad (3.5)$$

All cancellations occur inside the *Euler*[…] function. Its argument is:

$$\frac{zw}{Q} = \sum_{i=1}^{m} x_i \prod_{j=0}^{i-1} N_j \sum_{k=1}^{m} u_k \prod_{h=k+1}^{m+1} N_h \bigg/ \prod_{g=0}^{m+1} N_g \quad (3.6)$$

Each term, where $i > k$, becomes an integer and may be discarded. After removing them

Eq. (3.6) turns into:
$$\frac{zw}{Q} = \frac{x_m u_m}{N_m} + \sum_{i=1}^{m-1} \frac{x_i}{N_i}(u_i + \alpha_i) \quad (3.7)$$

Where
$$\alpha_i = \sum_{k=i+1}^{m} u_k \bigg/ \prod_{h=k+1}^{m+1} N_h \quad (3.8)$$

As has been expected the 1-D DFT turns into its M-D equivalent:

$$F_{1\text{-D}}(w) = F_{M\text{-D}}(u_1+\alpha_1, u_2+\alpha_2, ..., u_m) \quad (3.9)$$

## 4. Discussion

The main features of the new formulae are:

- **Index Reversal**. This was the key to separating the 1-D DFT into M dimensions. Thus, *e.g.* in 2-D, it means the following: If the image data are wrapped around the *Y* axis (*i.e.* stored by rows) then their single 1-D DFT should be wrapped around the *X* axis to produce the 2-D frequency data (stored by columns).

- **Limited Shifts**. For any number of dimensions, the α value is always below one. This basically means that this method's sampling points in the frequency domain lie within the pixel's size radius around the conventional (integer) sampling grid.

- **Proper Sampling**. The function remains properly sampled in the frequency domain since the sampling points stay regularly spaced. A minute stretching occurs along some axes, but, if our functions are over-sampled in the spatial domain, they should remain over-sampled on the stretched axes in the frequency domain.

- **Un-Orthogonality.** The single 1-D DFT of M-D spatial data *samples* the M-D frequency domain on an un-orthogonal (slightly skewed) grid, and vice-versa. That is: Applying the 1-D inverse DFT to the properly skewed M-D frequency data restores the spatial data. *Properly skewed* here means *the same as after applying the 1-D DFT*. In addition, the two above statements remain true if we swap the words *frequency* and *spatial*.

One of implications of this skewedness is that it has to be taken into account in such cases as *e.g.* transposition (columns into rows etc.). If we transpose our (supposedly orthogonal) *f(x,y)* and apply the 1-D DFT again, the overall picture of the frequency domain will look similar to the (transposed) original *F(u,v)*. However, because the *sampling* grid is skewed along the different axis, the actual pixel values differ (more or less) from the corresponding pixels in the original (before transposition) picture of frequency domain. The only un-skewable, and therefore unchangeable, frequency sample is the very first one (responsible for the average value across the all data).

## 5. Experiments

The 2-D formula (2.8) has been tested in a simple filtering experiment, the results of which are shown in the Figure 2. An 8-bit 256x256 colour image (Figure 2b) was used as a test picture. For clarity, it has only two colour components: red and green, to serve as the real and imaginary parts respectively. The image, containing some low (target-shaped background), mid (letters) and high (letter edges, the framing circle) frequency components, was subjected to a single 1-D FFT. This, according to the Eq. (2.8), produces the array of corresponding 2-D frequency data. The un-orthogonality of the array, for the purpose of displaying, is negligible (within a pixel) and was ignored (see Figure 2a. The picture data were transposed, and then centered by using the graphical device interface functions).

Figures 2c and 2d display the result of applying (convolving with) an ideal low-pass circular filter (2-D), and then restoring the filtered image by means of the single inverse 1-D FFT. This, as expected, has removed the high-frequency components: the image became blurred. Figures 2e and 2f show the case of ideal high-pass filtering (with exception of $F(0,0)$ to restore the average level), after which the background target has successfully disappeared.

Finally, here is the tally of number of multiplications required to calculate the Fourier transforms in the conventional way, and by this method of single 1-D FFT:

| 256 x 256 image | 2-D FFT | This method | Savings |
|---|---|---|---|
| Basic Operation (1-D FFT) | 5,116 | 2,359,292 | |
| Quantity | 512 | 1 | |
| Total Multiplications | 2,619,392 | 2,359,292 | 10% |

# Conclusion

- The restored images in Figure 2 have clearly confirmed the validity of Eq. (2.8).

- The use of single 1-D FFT instead of a 2-D FFT yields some minor advantages in the number of multiplications (*e.g*. 10% fewer for 256x256 images), plus significant simplification of required FFT routines.

- Also, the formulae give us some additional insights into the nature of these fabled multi-dimensional frequency spaces.

## *References*

1. R.C.Gonzales, and R.E.Woods. *Digital Image Processing*, Addison-Wesley, Reading, [1993].

2. J. C. Russ. *The Image Processing Handbook*, CRC Press, Boca Raton, [1995].

3. G. X. Ritter, and J. N. Wilson. *Handbook of Computer Vision Algorithms in Image Algebra*, CRC Press, Boca Raton, [1996].

4. A. R. Weeks. *Fundamentals of Electronic Image Processing*, SPIE / IEEE Press, Bellingham, [1996].

5. W. H. Press, S. A. Teukolsky, W. T. Vetterling, and B. P. Flannery. *Numerical Recipes in C*, Cambridge University Press, Cambridge, [1992].

**Andrew V. Batrac** received his degree in physics from the Moscow Engineering & Physics Institute, Russia, in 1986. He worked as research physicist at the Science & Research Centre on Technological Lasers, of the Russian Academy of Sciences. Currently he is with Polartechnics Ltd in Sydney, Australia, where he is involved in development of a skin-imaging device for the automated detection of melanoma.

Original state: 2-D Array of f(x, y).    Step 1. Re-Numbering as f(k), wrapping around Y.

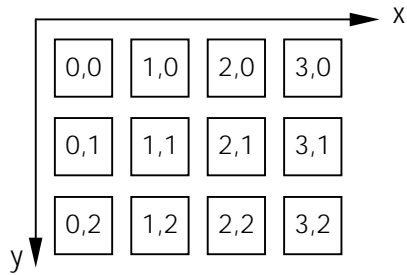   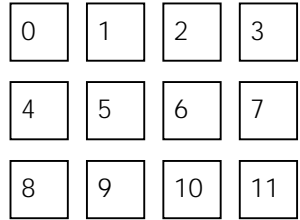

Step 1A. Re-Alignment (only visual, computer memory is already a line of data)

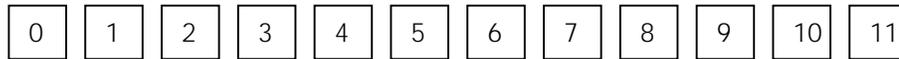

Step 2. Applying the 1-D DFT to f(k). The result: Array of F(j).

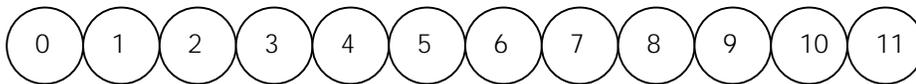

Step 2B. Arranging as 2-D array, wrapping around X (or U – analogue of X).

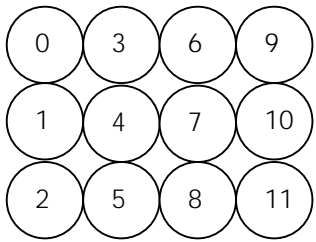

Step 3. Separation of x and y in the 1-D formula.
Result: Array of 2-D DFT, skewed along U axis.

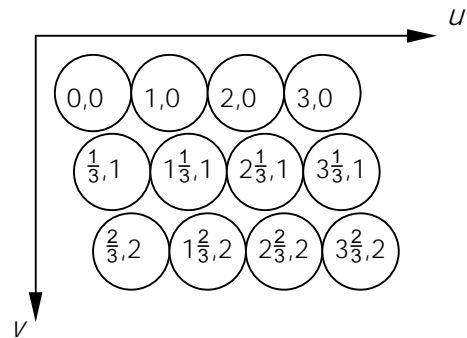

Figure 1.

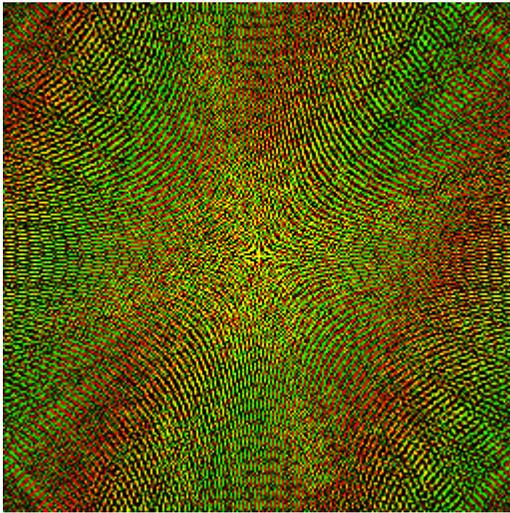
a) Single 1-D DFT of the image (b)

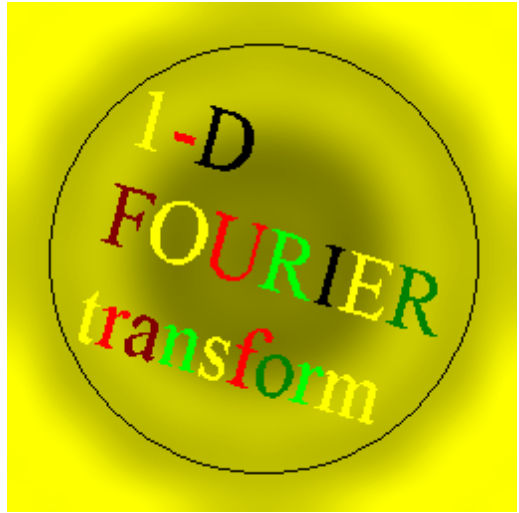
b) Original image 256x256.

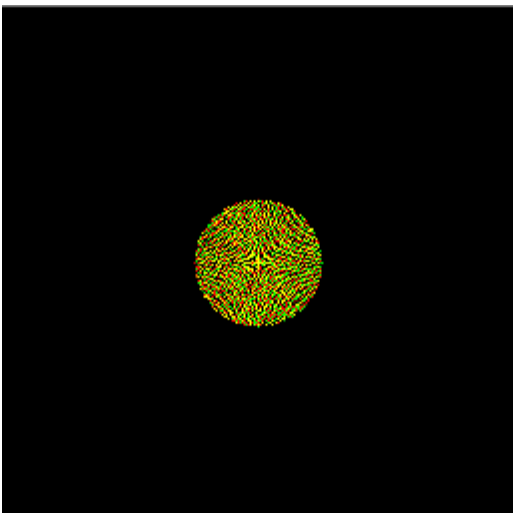
c) Applying an Ideal Low-Pass Filter

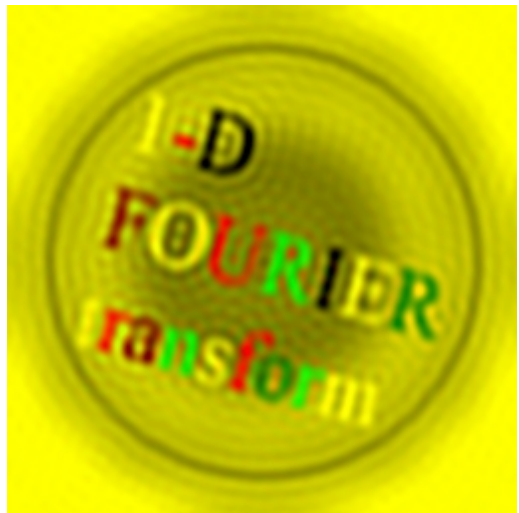
d) Restored image after ILPF

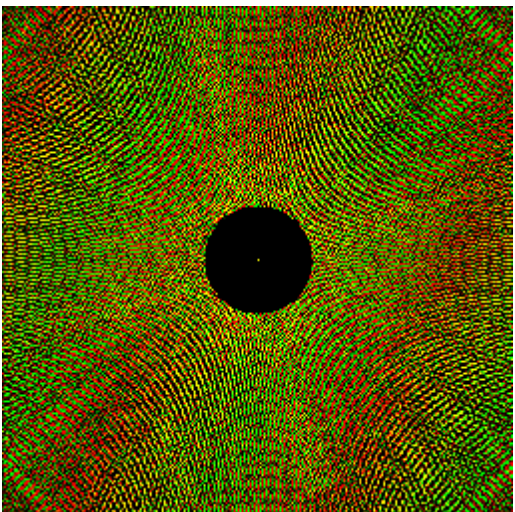
e) Applying an Ideal High-Pass Filter*.

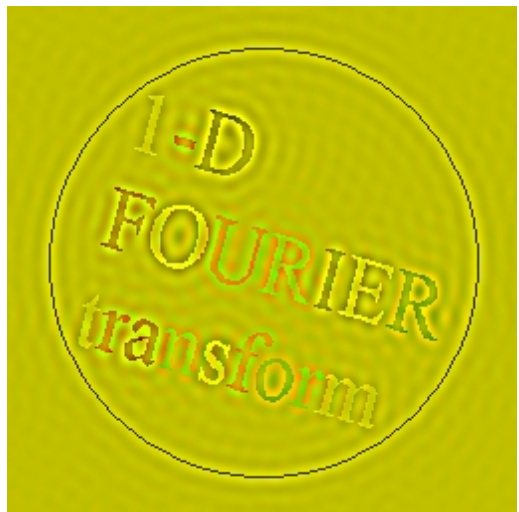
f) Restored image after IHPF*.

Figure 2.